\documentclass[reqno,12pt]{amsart}
\usepackage[mathscr]{eucal}
\usepackage{amsmath}
\usepackage{amssymb}
\usepackage{latexsym}
\usepackage{amsthm}
\usepackage{color}
\newtheorem{thm}{\indent Theorem}[section]
\newtheorem{cor}[thm]{\indent Corollary}
\newtheorem{lem}[thm]{\indent Lemma}

\newtheorem{dfn}{{\indent\bf Definition}}[section]

\newtheorem{expl}{{\indent\bf Example}}[section]

\newcommand{\lmx}{\left(\begin{matrix}}
\newcommand{\rmx}{\end{matrix}\right)}
\newcommand{\ldt}{\left|\begin{matrix}}
\newcommand{\rdt}{\end{matrix}\right|}

\newcommand{\td}{\tilde}

\newcommand{\nnm}{\nonumber}
\newcommand{\bbr}{{\mathbb R}}

\newcommand{\be}{\begin{equation}}
\newcommand{\ee}{\end{equation}}

\newcommand{\vfi}{\varphi}

\makeatletter
\numberwithin{equation}{section}
\begin{document}
\title [Calabi Hypersurfaces with parallel Fubini-Pick form]{{\large \bf Classification of Calabi Hypersurfaces with parallel Fubini-Pick form}}
\author [R.W. Xu and M.X. Lei] {Ruiwei Xu \quad\quad Miaoxin Lei }
\address{School of Mathematics and Information Sciences,
\newline \indent Henan Normal University, Xinxiang 453007, P. R. China
\newline \indent rwxu@htu.edu.cn;\; lmiaoxin@126.com}
\date{}
\footnotetext{$^1$\,The first author is partially supported by NSFC
 11871197 and 11671121.}
\maketitle
\renewcommand{\baselinestretch}{1}
\renewcommand{\arraystretch}{1.2}
\catcode`@=11 \@addtoreset{equation}{section} \catcode`@=12
\makeatletter
\renewcommand{\theequation}{\thesection.\arabic{equation}}
\@addtoreset{equation}{section} \makeatother
\noindent {\bf Abstract:} In this paper, we  present the classification of 2 and 3-dimensional Calabi hypersurfaces with parallel Fubini-Pick form with respect to the Levi-Civita connection of the Calabi metric.

\vskip 2pt\noindent {\bf
2000 AMS Classification:} Primary 53A15; Secondary 53C24, 53C42.
 \vskip 0.1pt

\section{Introduction}
In equiaffine differential geometry, the problem of classifying locally strong convex affine hypersurfaces with parallel Fubini-Pick form (also called cubic form) has been studied intensively, from the earlier beginning paper by Bokan-Nomizu-Simon \cite{BNS}, and then \cite{DV},\cite{DVY},\cite{HLSV}, to complete classification of Hu-Li-Vrancken \cite{HLV}. Here we recall results about the classification of locally strong convex affine hypersurfaces with parallel Fubini-Pick form $\nabla A=0$ with respect to the Levi-Civita connection of the affine Berwald-Blaschke metric. The condition $\nabla A=0$ implies that $M$ is an affine hypersphere with constant affine scalar curvature. Thus Theorem 1 of Li and Penn \cite{LP} (see also Theorem 3.7 in \cite{LSZH}) can be restated as follows:

\begin{thm}{\cite{LP}}
Let $x:M\rightarrow A^3$ be a locally strongly convex affine surface with $\nabla A=0$. Then, up to an affine transformation, $x(M)$ lies on the surface $x_1x_2x_3=1$ or strongly convex quadric.
\end{thm}

\noindent The classification of 3-dimensional affine hypersurfaces with parallel Fubini-Pick form, due to Dillen and Vrancken \cite{DV}.

\begin{thm}{\cite{DV}}
Let $x:M\rightarrow A^4$ be a locally strongly convex affine hypersurface with $\nabla A=0$. Then, up to an affine transformation, either $x(M)$ is an open part of a locally strongly convex quadric or $x(M)$ is an open part of one of the following two hypersurfaces:\\
 (i)\; $x_1x_2x_3x_4=1$;\\
 (ii)\;$(x_1^2-x_2^2-x_3^2)^3x_4^2=1$.
\end{thm}
\noindent In \cite{HLV}, Hu-Li-Vrancken introduced some typical examples and gave the complete classification of locally strongly convex affine hypersurfaces of $R^{n+1}$ with parallel cubic form with respect to the Levi-Civita connection of the affine Berwald-Blaschke metric.

In centroaffine differential geometry, Cheng-Hu-Moruz \cite{CHM} obtained a complete classification of locally strong convex centroaffine hypersurfaces with parallel cubic form. On the other hand, Liu and Wang \cite{LW1} gave the classification of the centroaffine surfaces with parallel traceless cubic form relative to the Levi-Civita connection. In \cite{CH}, Cheng and Hu established a general inequality for locally strongly convex centroaffine hypersurfaces in $R^{n+1}$ involving the norm of the covariant differentiation of both the difference tensor and the Tchebychev vector field. Applying the classification result of \cite{CHM}, Cheng-Hu \cite{CH} completely classified locally strongly convex centroaffine hypersurfaces with parallel traceless difference tensor.

A centroaffine hypersurface is said to be \emph{Canonical} if its Blaschke metric is flat and its Fubini-Pick form is parallel with respect to its Blaschke metric. In \cite{LW}, Li and Wang classified the Canonical centroaffine hypersurfaces in $R^{n+1}$. In this paper, we first classify the canonical Calabi hypersurfaces in Calabi geometry. As a corollary, we classify Calabi surfaces with parallel Fubini-Pick form with respect to the Levi-Civita connection of the Calabi metric.

\begin{thm}  Let $f$ be a smooth strictly convex function on a domain $\Omega\in R^n$. If  its graph $M=\{(x,f(x)) \;| \;   x\in\Omega\}$ has a flat Calabi metric and parallel cubic form.
Then $M$ is Calabi affine equivalent to an open part of the following hypersurfaces:
\\ (i) elliptic paraboloid; or
\\ (ii) the hypersurfaces $Q(c_1,\cdots,c_r;n),\, 1\leq r \leq n.$
\end{thm}
The definitions of \emph{Calabi affine equivalent} and \emph{hypersurfaces $Q(c_1,\cdots,c_r;n)$} will be given in Section 2 (see Definition 2.1 and Example 2.1, respectively).

\begin{cor}
 Let $f$ be a smooth strictly convex function on a domain $\Omega\in R^2$. If  its graph $M=\{(x,f(x)) \;| \;   x\in\Omega\}$  has parallel cubic form, $M$ is Calabi affine equivalent to an open part of the following surfaces:
\\ (i) elliptic paraboloid; or
\\ (ii) the surfaces $Q(c_1,\cdots,c_r;2),\, 1\leq r \leq 2.$
\end{cor}

Motivated by above classification results in equiaffine differential geometry and centroaffine differential geometry, we present the classification of 3-dimensional Calabi hypersurfaces with parallel Fubini-Pick form with respect to the Levi-Civita connection of the
Calabi metric in Section 4 and Section 5.

\begin{thm} \label{main2} Let $f$ be a smooth strictly convex function on a domain $\Omega\in R^3$. If  its graph $M=\{(x,f(x)) \;| \;   x\in\Omega\}$  has  parallel cubic form.
Then $M$ is Calabi affine equivalent to an open part of one of the following three types of hypersurfaces:
\\ (i) elliptic paraboloid; or
\\ (ii) the hypersurfaces
   $Q(c_1,\cdots,c_r;3),\, 1\leq r \leq 3;$ or
\\ (iii) the hypersurface
$$x_4= -\frac{1}{2c^2}\ln(x_1^2-(x_2^2+x_3^2)),$$
 where the constant $-2c^2$ is the scalar curvature of $M$.
\end{thm}

\textbf{Remark:} (1) In \cite{XL}, Xu-Li proved that

\begin{thm}\cite{XL}
Let $M^n (n\geq3)$ be a Calabi complete Tchebychev affine K\"ahler hypersurface with nonnegative Ricci curvature. Then it must be Calabi affine equivalent to either an elliptic paraboloid or one of the hypersurfaces $Q(c_1,...,c_r;n)$.
\end{thm}
Case (iii) of Theorem 1.5 shows that there exists a class of Calabi complete Tchebychev affine K\"ahler hypersurfaces with negative Ricci curvature. Thus the restriction on Ricci curvature in Theorem 1.6 is essential.

(2) The Euler-Lagrange equation of the volume variation with respect to the Calabi metric can be written as the following fourth order PDE (see \cite{L} or \cite{LXSJ})
\be\label{aeh} \Delta \ln \det\left(\frac{\partial^2f}{\partial x_i\partial x_j}\right)=0,\ee
where $\Delta$ is the Laplacian of the  Calabi metric $G=\sum \frac{\partial^2f}{\partial x_i\partial x_j}dx_idx_j$. Graph hypersurfaces $M^n$ defined by solutions of \eqref{aeh} are called \emph{affine extremal hypersurfaces}.
Case (iii) of Theorem 1.5 shows that there exists a class of new Euclidean complete and Calabi complete affine extremal hypersurfaces.

\section{Preliminaries}
\subsection {Calabi geometry.}

In this section, we shall show some basic facts for the Calabi geometry, see \cite{Ca} or \cite{P}. Let $f$ be a strictly convex $C^\infty$-function on a domain $\Omega \in \mathbb{R}^n.$ Consider the graph hypersurface
\be\label{2.1} M^n:=\{(x_i,f(x_i))\;|\;x_{n+1}=f(x_1,\cdots,x_n), (x_1,\cdots,x_n)\in\Omega\}.\ee
This Calabi geometry can also be essentially realised as a very special relative affine geometry of the graph hypersurface (\ref{2.1})
by choosing the so-called {\em Calabi affine normalization} with $$Y =(0,0,\cdots,1)^t\in R^{n+1}$$
being the fixed relative affine normal vector field which we call the {\em Calabi affine normal}.

For the position vector $x=(x_1,\cdots,x_n,f(x_1,\cdots,x_n))$ we have the decomposition
\be x_{ij}=c_{ij}^k
x_*\frac{\partial }{\partial x_k}+f_{ij}Y,\ee
with respect to the bundle decomposition $\bbr^{n+1}=x_*TM^n\oplus\bbr\cdot Y$, where the induced affine connection $c^k_{ij}\equiv 0$. It follows that the relative affine metric is nothing but the Calabi metric
$$ G=\sum \frac{\partial^2f}{\partial x_i\partial x_j}dx_idx_j.$$
The Levi-Civita connection with respect to the metric $G$ has the Christoffel symbols
$$\Gamma_{ij}^k=\frac{1}{2}\sum f^{kl}f_{ijl},$$
where and hereafter
$$f_{ijk}=\frac{\partial^3f}{\partial x_i\partial x_j\partial x_k}, \quad (f^{kl})=(f_{ij})^{-1}. $$
Then we can rewrite the Gauss structure equation as follows:\be x_{,ij}=\sum A_{ij}^kx_*\frac{\partial }{\partial x_k}+G_{ij}Y.\ee
The  Fubini-Pick tensor (also called cubic form) $A_{ijk}$ and the Weingarten tensor satisfy
\be A_{ijk}=A_{ij}^lG_{kl}=-\frac{1}{2}f_{ijk},\quad B_{ij}=0,\ee
which means $A_{ijk}$ are symmetric in all indexes.
Classically, the tangent vector field
\be T:=\frac{1}{n}\sum G^{kl}G^{ij}A_{ijk}\frac{\partial}{\partial x_l}=-\frac{1}{2n}\sum f^{kl}f^{ij}f_{ijk}\frac{\partial}{\partial x_l}\ee
is called the  \emph{Tchebychev vector field} of the hypersurface $M^n$, and the invariant function
\be J:=\frac{1}{n(n-1)}\sum G^{il}G^{jp}G^{kq}A_{ijk}A_{lpq}=\frac{1}{4n(n-1)}\sum f^{il}f^{jp}f^{kq}f_{ijk}f_{lpq}\ee
is named as the  \emph{relative Pick invariant}  of $M^n$. As in the report \cite{Li},  $M^n$ is called a \emph{Tchebychev affine K\"{a}hler hypersurface},
if the Tchebychev vector field $T$ is parallel with respect to the Calabi metric $G$.
\\The \emph{Gauss integrability conditions}  and the  \emph{Codazzi  equations}  read
\be \label{2.7} R_{ijkl}=\sum f^{mh}(A_{jkm}A_{hil}-A_{ikm}A_{hjl}),\ee
\be \label{2.8} A_{ijk,l}=A_{ijl,k}.\ee
From (\ref{2.7}) we get the  \emph{Ricci tensor}
\be\label{2.9} R_{ik}=\sum f^{jl}f^{mh}(A_{jkm}A_{hil}-A_{ikm}A_{hjl}).\ee
Thus, the scalar curvature is given by
\be\label{2.10} R=n(n-1)J-n^2|T|^2.\ee

Using Ricci identity, (\ref{2.7}) and (\ref{2.8}), we have two useful formulas.

\begin{lem}
For a Calabi hypersurface, the following formulas hold
\be\label{2.11}\frac{1}{2}\Delta|T|^2=\sum T_{i,j}^2+\sum T_{i}T_{j,ji}+\sum R_{ij}T_{i}T_{j},\ee
\be\label{2.12} \frac{n(n-1)}{2}\Delta J
     =\sum(A_{ijk,l})^2+\sum A_{ijk}A_{lli,jk}+\sum(R_{ijkl})^2+\sum R_{ij}A_{ipq}A_{jpq}.\ee
\end{lem}

Let $A(n+1)$ be the group of $(n+1)$-dimensional affine transformations on $R^{n+1}$. Then $A(n+1)=GL(n+1)\ltimes R^{n+1}$, the semi-direct product of the general linear group $GL(n+1)$ and the group $R^{n+1}$ of all the parallel transports on $R^{n+1}$. Define
\be
SA(n+1)=\{\phi=(M,b)\in A(n+1)=GL(n+1)\ltimes R^{n+1};\ M(Y)=Y\}
\ee
where $Y=(0,\cdots,0,1)^t$ is the Calabi affine normal. Then the subgroup $SA$ consists of all the transformations $\phi$ of the following type:
\begin{align}
X:&=(X^j,X^{n+1})^t\equiv(X^1,\cdots, X^n,X^{n+1})^t\nnm\\
&\mapsto \phi(X):=\lmx a^i_j&0\\a^{n+1}_j&1\rmx X+b,\quad \forall\, X\in R^{n+1},
\end{align}
for some $(a^i_j)\in A(n)$, constants $a^{n+1}_j$ ($j=1,\cdots,n$) and some constant vector $b\in R^{n+1}$. Clearly, the Calabi metric $G$ is invariant under the action of $SA(n+1)$ on the graph hypersurfaces or, equivalently, under the induced action of $SA(n+1)$ on the strictly convex functions, which is naturally defined to be the composition of the following maps:
$$
f\mapsto (x_i,f(x_i))\mapsto(\td x_i,\td f(\td x_i)):=\phi(x_i,f(x_i))\mapsto \phi(f):=\td f,\quad\forall\,\phi\in SA(n+1).
$$

\begin{dfn}\cite{XL}
Two graph hypersurfaces $(x_i,f(x_i))$ and $(\td x_i,\td f(\td x_i))$, defined respectively in domains $\Omega,\td\Omega\subset R^n$, are called Calabi-affine equivalent if they differ only by an affine transformation $\phi\in SA(n+1)$.
\end{dfn}

Accordingly, we have

\begin{dfn}\cite{XL}
Two smooth functions $f$ and $\td f$ respectively defined on domains $\Omega,\td\Omega\subset R^n$ are called affine equivalent $($related with an affine transformation $\vfi\in A(n))$ if there exist some constants  $a^{n+1}_1,\cdots,a^{n+1}_n,b^{n+1}\in R$ such that $\vfi(\Omega)=\td\Omega$ and
\be\td f(\vfi(x_j)) =f(x_j)+\sum a^{n+1}_jx_j+b^{n+1}\quad \text{for all}\ (x_j)\in\Omega.\ee
\end{dfn}

Clearly, the above two definitions are equivalent to each other.
\subsection {Canonical Calabi hypersurfaces}
A Calabi hypersurface is called \emph{canonical} if its Fubini-Pick form is parallel with respect to the Levi-Civita connection of the Calabi metric and its Calabi metric $G$ is flat. In \cite{XL} Xu and Li introduced a large class of new canonical Calabi hypersurfaces, which are denoted by $Q(c_1,\cdots,c_r;n)$. It turns out that these new examples are all Euclidean complete and Calabi complete.

\begin{expl} \cite{XL}  Given the dimension $n$ and let $1\leq r\leq n$. For any positive numbers $c_1,\cdots, c_r$, define
$$\Omega_{c_1,\cdots,c_r;n}=\{(x_1,\cdots,x_n);\ x_1>0,\cdots,x_r>0\}$$
and consider the following smooth functions
\begin{align}
f(x_1,\cdots,x_n)\equiv&\; Q(c_1,c_2,\cdots,c_r;n)(x_1,\cdots,x_n)\nnm\\
:=&-\sum_{i=1}^r c_i\ln x_i+\sum_{j=r+1}^n \frac12 x_j^2,\qquad (x_1,\cdots,x_n)\in\Omega_{c_1,\cdots,c_r;n}.\label{2.14}
\end{align}
\end{expl}

\section{Proof of the Theorem 1.3 and Corollary 1.4}

\begin{lem}\label{lemma3.1}
Let $ \{A^i\}_{1\leq i\leq n}$ be real symmetric matrices satisfying $A^iA^j=A^jA^i$, $\forall{1\leq i,j\leq n}$. Then there exists an orthogonal matrix such that matrices $ \{A^i\}_{1\leq i\leq n} $ can be simultaneously diagonalized.
\end{lem}

{\bf Proof.}
As we know the conclusion obviously holds for the case of $n=2$. Now we assume that it holds for $n=k$. Namely, there is an orthogonal matrix $P$ such that $A^1,\cdots,A^k$ can be simultaneously diagonalized:
\be\label{3.1} PA^iP^{-1}= diag(\lambda^i_1E_{n^i_1},\cdots,\lambda^i_sE_{n^i_s}),\quad 1\leq i\leq k,\ee
where $\lambda^i_1,\cdots,\lambda^i_s$ are different eigenvalues for any fixed $i$.

In the next we will prove the conclusion still holds for the case of $n=k+1$.
Since $A^iA^{k+1}=A^{k+1}A^i,\forall\,1\leq i\leq k,$ we obtain \be\label{3.2} (PA^iP^{-1})(PA^{k+1}P^{-1})=(PA^{k+1}P^{-1})
(PA^iP^{-1}).\ee
Denote $$B^{k+1}:=PA^{k+1}P^{-1},$$
where $B^{k+1}$ is a real symmetric matrix.
From (\ref{3.1}) and (\ref{3.2}) we have
\be\label{3.3}\lambda^i_pB^{k+1}_{pq}=B^{k+1}_{pq}
\lambda^i_q,\quad\forall p,q,\quad 1\leq i\leq k.\ee

Fix an arbitrary index $i\in\{1,\cdots,k\}$. If $\lambda^i_p\neq \lambda^i_q$ for some indices $p\neq q$ then, by (\ref{3.3}), it must hold that $B^{k+1}_{pq}=0$. Therefore, for any pair of indices $p\neq q$, if $B^{k+1}_{pq}\neq 0$, then it holds that $\lambda^i_{p}=\lambda^i_q$ for each $i=1,\cdots,k$.
Thus we get
$$B^{k+1}=diag(B^{k+1}_{n^{k+1}_1},\cdots,B^{k+1}_{n^{k+1}_r}),$$
where $B^{k+1}_{n^{k+1}_j}(1\leq j\leq r)$ are real symmetric matrices of order $n^{k+1}_j$, and, for any fixed $1\leq j\leq r$ and $1\leq i\leq k$, the $n^{k+1}_1+\cdots+n^{k+1}_{j-1}+1$ th to $n^{k+1}_1+\cdots+n^{k+1}_{j}$ th eigenvalues of $A^i$ are equal. Thus there are a set of orthogonal matrices $R_{n^{k+1}_j},1\leq j\leq r$ such that
$R_{n^{k+1}_j}B^{k+1}_{n^{k+1}_j}R^{-1}_{n^{k+1}_j}
$ are diagonal matrices.
Let $$R=diag(R_{n^{k+1}_1},\cdots,R_{n^{k+1}_r}).$$
Then the real symmetric matrices $A^1,\cdots,A^{k+1}$ can be simultaneously diagonalized by orthogonal matrix $RP$. \hfill$\Box$

\textbf{Proof of the Theorem 1.3.}
The canonical Calabi hypersurface means
$$ \nabla A=0 \qquad \text{and} \qquad  R_{ijkl}=0.$$
Hence $M^n$ locally is a Euclidean space. We choose local coordinates $\{u^1,\cdots,u^n\}$ such that the Calabi metric is given by $G=\sum (du^i)^2$, and
$A_{ijk}=const$ in this coordinates.  We consider the following two subcases:

{\bf Case 1.} \; $A_{ijk}=0,\forall i,j,k.$ Obviously, in this case, $M^n$ is an open part of elliptic paraboloid.

{\bf Case 2.} \; Otherwise. Let $p\in M^n$ be a fixed point with coordinates $(0,\cdots,0)$. Choose the local orthonormal frame field $ e_i=\frac{\partial}{\partial u^i} $,  and $e=(0,\cdots,0,1)$ on $M^n$. Let $\{\omega ^i\}$ be the dual frame field of $\{e_i\}$.
Denote $$A^{(k)}:=A^{e_k}=\sum A_{ij}^{(k)}du^idu^j,$$
$$ A_{ij}^{(k)}:=A(e_i,e_j,e_k)\equiv A_{ijk}.$$
By the Gauss integrability conditions (\ref{2.7}) and the flatness of the metric $G$, we have
\be\sum A_{iml}A_{jmk}-\sum A_{imk}A_{jml}=0,\ee
which means the following matrix equalities:
$$(A_{ij}^{(k)})(A_{ij}^{(l)})=(A_{ij}^{(l)})(A_{ij}^{(k)}), \qquad \forall 1\leq k,l\leq n.$$
By Lemma 3.1, we get that matrices $(A_{ij}^{(k)})$ can be simultaneously diagonalized. There exists an orthogonal constant matrix $C=(c_{ij}),$ for any fixed $1\leq k\leq n,$ such that
$$(A_{\bar{i} \bar{j}}^{(k)})
=C(A_{ij}^{(k)})C^{-1}= diag(\lambda_1^k,\lambda_2^k,\cdots,\lambda_n^k).$$
Here $A_{\bar {i}\bar {j}}^{(k)}=A({\bar {e}_i},{\bar {e}_j},e_k)$
and $\bar{e}_i=\sum c_{ij}e_j,1\leq i\leq n,$  then
$$\bar A_{ijk}:=A({\bar {e}_i},{\bar {e}_j},{\bar {e}_k})=\sum c_{kl}A(\bar{e}_i,\bar{e}_j,e_l)=\sum c_{kl} \lambda _i^l \delta _{ij}.$$
Since the matrices $(\bar{A}_{ijk})$ are symmetric in all indexes, we get:
\begin{align}
\bar{A}_{ijk}=
\begin{cases}
\bar{A}_{iii},& 1\leq i=j=k\leq n,\\
 0,& otherwise.\tag{3.5}
\end{cases}
\end{align}
\\From $dx=\omega^ie_i=\bar{\omega}^i\bar{e}_i,$ we can get $\bar{\omega}^i=\sum c^{ij}\omega^j=\sum c^{ij}du^j,$ where $(c^{ij})$ denotes the inverse matrix of $(c_{ij}).$ Let $\bar{u}^i=\sum c_{ij}u^j, 1\leq i\leq n,$ then $(\bar{u}^1,\cdots,\bar{u}^n)$ are new Euclidean coordinates of $M^n$, such that $\frac{\partial}{\partial \bar{u}^i}=\bar{e}_i,1\leq i\leq n.$
Under these new coordinates, the tensor $\bar{A}$ is expressed as (3.5), thus we have:
\begin{align}
\begin{cases}
d\bar{e}_i=\sum \bar{\omega}_i^j\bar{e}_j+d\bar{u}^ie, \qquad 1\leq i\leq n ,\\
 dx=\sum d\bar{u}^i \bar{e}_i.\tag{3.6}
\end{cases}
\end{align}
Since the Calabi metric is flat and $\frac{\partial}{\partial \bar{u}^i}$ are orthonormal, we obtain $\bar{\omega}_i^j=\bar{A}_{ijk}d\bar{u}^k.$
Assume that $x\in M^n$ is an arbitrary point with coordinates $(v^1,\cdots,v^n).$ We draw a curve connecting $p$ and $x$
$$\bar{u}^i(t)=v^it,\qquad 0\leq t\leq 1.$$
Along this curve the equations (3.6) become
\begin{align}
\begin{cases}
\frac{d\bar{e}_i}{dt}=\bar{A}_{iii}v^i\bar{e}_i+v^ie, \qquad 1\leq i\leq n ,\\
\frac{dx}{dt}=\sum v^i \bar{e}_i.\tag{3.7}
\end{cases}
\end{align}
Consider the ordinary differential equation
$$\frac{du}{dt}=au+b.$$
It is easy to find out its solution
\begin{align}
 u(t)=
\begin{cases}
(u(0)+\frac{b}{a})e^{at}-\frac{b}{a}, &\qquad a\neq0,\\
 u(0)+bt,&\qquad a=0.
\end{cases}\notag
\end{align}
We may assume that $\bar{e}_i(0)=(0,\cdots,1,\cdots,0),1\leq i\leq n,$ where 1 is on i-th entry and $\bar{A}_{iii}\geq0$ at point $p$. By an arrangement, we can get
\begin{align}
\begin{cases}
\bar{A}_{iii}>0,&\qquad 1\leq i\leq r;\\
\bar{A}_{jjj}=0,&\qquad r+1\leq j\leq n,\tag{3.8}
\end{cases}
\end{align}
where $1\leq r\leq n$ and $r=n$ means that $\bar{A}_{iii}>0$ for all $1\leq i\leq n.$
\\Without loss of generality, we assume $v^i>0,1\leq i\leq n.$ Solve equations (3.7), we obtain:
\begin{align}
\begin{cases}
\bar{e}_{i}=\exp(\bar{A}_{iii}v^it)\bar{e}_i(0)+
\frac{1}{\bar{A}_{iii}}\exp(\bar{A}_{iii}v^it)e
-\frac{1}{\bar{A}_{iii}}e,&\qquad 1\leq i\leq r;\\
\bar{e}_j=\bar{e}_j(0)+v^jte,&\qquad r+1\leq j\leq n;\\
x(t)=x(0)+\int_0^tv^i\bar{e}_i(s)ds.\tag{3.9}
\end{cases}
\end{align}
Thus
\begin{align}
x(t)=&\;x(0)+\sum_{i=1}^r \frac{1}{\bar{A}_{iii}}[\exp(\bar{A}_{iii}v^it)-1]\bar{e}_i(0)
+\sum_{i=1}^r \frac{1}{\bar{A}^2_{iii}}[\exp(\bar{A}_{iii}v^it)-1]e\notag\\
&-\sum_{i=1}^r \frac{1}{\bar{A}_{iii}}v^ite
+\sum_{j=r+1}^n v^jt\bar{e}_j(0)+\sum_{j=r+1}^n\frac{1}{2}(v^j)^2t^2e.\tag{3.10}
\end{align}
Evaluate (3.10) at $t=1$ we have:
\begin{align*}
x_i=&\;x_i(0)+\frac{1}{\bar{A}_{iii}}[\exp(\bar{A}_{iii}v^i)-1],\qquad \qquad 1\leq i\leq r;\\
x_j=&\;x_j(0)+v^j, \qquad\qquad \qquad\qquad \qquad\qquad r+1\leq j\leq n;\\
x_{n+1}=&\;x_{n+1}(0)+\sum_{i=1}^r\left(\frac{1}{\bar{A}^2_{iii}}
[\exp(\bar{A}_{iii}v^i)-1] - \frac{1}{\bar{A}_{iii}}v^i\right)
+\sum_{j=r+1}^n\frac{1}{2}(v^j)^2. \tag{3.11}
\end{align*}
Inserting $x_i$ and $x_j$ into $x_{n+1}$, we find
$$x_{n+1}=\sum_{i=1}^r\frac{1}{\bar{A}_{iii}}x_i-
\sum_{i=1}^r\frac{1}{\bar{A}^2_{iii}}\ln(\bar{A}_{iii}x_i+1)
+\frac{1}{2}\sum_{j=r+1}^n(x_j)^2.$$
It is easy to find that $x_{n+1}$ is affine equivalent to
$$x_{n+1}=-\sum_{i=1}^r\frac{1}{\bar{A}^2_{iii}}\ln x_i+\sum_{j=r+1}^n\frac{1}{2}(x_j)^2.\eqno(3.12)$$
This completes the proof of theorem 1.3. \hfill $\Box$
\vskip 0.1in {\bf Proof of Corollary 1.4}

By $\nabla A=0$, the definition of the Tchebychev vector field $T$ and the Pick invariant $J$, we can get:
$$\nabla T=0 \qquad \text{and} \qquad J=const.$$
It follows that $|T|=const.$

{\bf Case 1.} \; $|T|=0$. It means that
$$\det(f_{ij})=const>0,\eqno(3.13)$$
and the Ricci formula
$$R_{ij}=A_{iml}A_{jml}-A_{ijm}A_{mll}=A_{iml}A_{jml}.\eqno(3.14)$$
Hence, by (3.14), (\ref{2.12}) and $\nabla A=0$, we have
\begin{align}
\frac{n(n-1)}{2}\Delta J                        =\sum(R_{ij})^2+\sum(R_{ijkl})^2.\tag{3.15}
\end{align}
It follows that $R_{ijkl}=0$. Then, by (\ref{2.10}), we obtain the relative Pick invariant
$$n(n-1)J=R+n^2|T|^2=0.$$
Thus $f$ is a strictly convex quadratic function.

{\bf Case 2.} \;$|T|=const>0.$ In this case, we can choose an orthonormal frame field $\{\tilde{e}_1,\tilde{e}_2\}$ on $M^2$ with $\tilde{e}_1=\frac{T}{|T|},$ where $\nabla\tilde{e}_1=0,$ since $\nabla T=0.$ From the definition of the Riemannian curvature tensor, we get£º$$R_{ijkl}=0.\eqno(3.16)$$
Thus, by Theorem 1.3, we complete
the proof of the corollary 1.4.\hfill $\Box$
\section{The classification of 3-dimension case}
\subsection{Elementary discussions in terms of a typical basis}
Now, we fix a point $p \in M^n$. For subsequent purpose, we will review the well known construction of a typical orthonormal basis for $T_p M^n$, which was introduced by Ejiri and has been widely applied, and proved to be very useful for various situations, see e.g., \cite{HLSV}, \cite{LV} and \cite{CHM}. The idea is to
construct from the $(1, 2)$ tensor $A$ a self adjoint operator at a point; then one extends the eigenbasis to a local field.
Let $p \in M^n$ and $U_p M^n = \{v \in T_pM^n \;|\; G(v, v) = 1\}$. Since $M^n$ is locally strong convex, $U_pM^n$ is compact. We define a function $F$ on $U_pM^n$ by $F (v) = A(v,v,v)$.
Then there is an element $e_1\in U_pM^n$ at which the function $F(v)$ attains an absolute maximum, denoted by  $\mu_1$. Then we have the following lemma. For its proof, we refer the reader to \cite{HLSV} or \cite{LSZH}.

\begin{lem}\label{lemma4.1}  There exists an orthonormal basis $\{e_1,\cdots,e_n\}$ of $T_pM^n$ such that the following hold: \\
(i) $A(e_1,e_i,e_j)=\mu_i\delta_{ij}$, for $i=1,\cdots,n$. \\
(ii) $\mu_1\geq2\mu_i,$ for $i\geq2$. If $\mu_1=2\mu_i$, then $A(e_i,e_i,e_i) =0$.
\end{lem}
Consider the function
$$F(v)=A(v,v,v)\qquad \;\text{on}\;\;\; U_p M^n.$$
Let $e_1\in U_pM^n$ be a vector at which $F(v)$ attains an absolute maximum $A_{111}(\geq0)$. From Lemma 4.1, we can further choose $e_2,\cdots, e_n$ such that $\{e_1,\cdots,e_n\}$ form an orthonormal basis of $T_pM^n,$ which possesses the following properties:
\begin{align}
&G(e_i,e_j)=\delta_{ij},\qquad A_{1ij}=\mu _i\delta_{ij},\qquad 1\leq i,j\leq n;\notag\\
&\mu_1\geq2\mu_i \text { and if }  \mu_1 =2\mu_i,\;  \text {then} \; A(e_i,e_i,e_i)=0 \; \text {for} \;i\geq 2.\notag
\end{align}
Using $\nabla A=0$ and the Ricci identity, for $i\geq2$, we have
\begin{align}\label{4.1}
0&=A_{11i,1i}-A_{11i,i1}=2A_{p1i}R_{p11i}+A_{11p}R_{pi1i}\notag\\
 &=\mu_i(\mu_1-2\mu_i)(\mu_i-\mu_1).
\end{align}

Therefore we have the following lemma.

\begin{lem}\label{lemma4.2} Let $M^n$ be a Calabi hypersurface with parallel Fubini-Pick form. Then, for every point $p\in M^n$, there exists an orthonormal basis $\{e_j\}_{1\leq j\leq n}$ of $T_pM^n$ (if necessary, we rearrange the order), satisfying $A(e_1,e_j)=\mu_je_{j}$, and there exists a number $i$, $0\leq i\leq n$, such that
$$\mu_2=\mu_3=\cdots=\mu_i=\frac{1}{2}\mu_1;\quad \mu_{i+1}=\cdots=\mu_n=0.$$
\end{lem}
Therefore, for a strictly convex Calabi hypersurface with parallel Fubini-Pick form, we have to deal with $(n+1)$ cases as follows: \\

\noindent\textbf{Case} $\mathfrak{C}_0.$  $\mu_1=0$.

\noindent\textbf{Case} $\mathfrak{C}_1.$ $\mu_1>0;\mu_2=\mu_3=\cdots=\mu_n=0.$

\noindent\textbf{Case }$\mathfrak{C}_i. $ $\mu_2=\mu_3=\cdots=\mu_i=\frac{1}{2}\mu_1>0;\quad \mu_{i+1}=\cdots=\mu_n=0 \quad \text{for}\;\; 2\leq i\leq n-1.$

\noindent\textbf{Case} $\mathfrak{C}_n. $ $\mu_2=\mu_3=\cdots=\mu_n= \frac{1}{2}\mu_1>0.$

When working at the point $p\in M^n$, we will always assume that an orthonormal basis is chosen such that Lemma 4.1 is satisfied.

\subsection{The settlement of the Cases $\mathfrak{C}_0$ and $\mathfrak{C}_n$}
Firstly, about the Case $\mathfrak{C}_0$, we have the following lemma.
\begin{lem}\label{lemma4.3}   If the Case $\mathfrak{C}_0$  occurs, then $M^n$ is an open part of elliptic paraboloid.
\end{lem}
\textbf{ Proof.}  If $\mu_1=0$, then \be \label{4.2} A(v,v,v)=0\quad  \text{for any} \quad v\in U_pM^n.\ee
Put $v=\frac{1}{\sqrt 2} (e_i+e_j)\in U_pM^n$ in (\ref{4.2}), then
$$0=A(e_i,e_i,e_j)+A(e_i,e_j,e_j).$$
On the other hand, put $v=\frac{1}{\sqrt 2} (e_i-e_j)\in U_pM^n$ in (\ref{4.2}), then
$$0=-A(e_i,e_i,e_j)+A(e_i,e_j,e_j).$$
Thus we have
$$A(e_i,e_i,e_j)=0,\quad\forall i,j.$$
From $0=\frac{1}{2} A(e_i+e_k,e_i+e_k,e_j)$, we have
$$A(e_i,e_j,e_k)=0,\quad\forall i,j,k. $$
Therefore $J\equiv0$, and $M^n$ is an open part of elliptic paraboloid.\hfill $\Box$

Secondly, we have the following important observation:
\begin{lem}\label{lemma4.4}   The Case $\mathfrak{C}_n$ does not occur.
\end{lem}
\textbf{Proof.}  Assume that this case does occur.
For any $i\geq2$, $\mu_i=\frac{1}{2}\mu_1>0$, then $A(e_1,v,v)=\frac{1}{2}\mu_1$ and $A(v,v,v)=0$ for any $v\in \{e_1^\perp\}\bigcap U_pM^n$. From the proof of Lemma 4.3, we see that
$$A(e_i,e_j,e_k)=0, \quad 2\leq i,j,k\leq n.$$
Then, for any unit vector $v\in \{e_1^\perp\}\bigcap U_pM^n$, we have
\be\label{4.3} A(e_1,e_1)=\mu_1 e_1,\quad A(e_1,v)=\frac{1}{2}\mu_1 v,\quad A(v,v)=\frac{1}{2}\mu_1 e_1. \ee
By $\nabla A=0$, we know that the curvature operator of Levi-Civita $R$  and  Fubini-Pick tensor $A$ satisfy
\be\label{4.4} R(e_1,v)A(v,v)=2A(R(e_1,v)v,v).\ee
By (\ref{4.4}), (\ref{4.3}) and (\ref{2.7}), we get $\mu_1=0$. This contradiction completes the proof of Lemma 4.4.\hfill $\Box$

In the following we only consider 3-dimensional Calabi hypersurfaces with parallel Fubini-Pick form.
Therefore, we only need to deal with the
\textbf{Case }$\mathfrak{C}_1$ and \textbf{Case }$\mathfrak{C}_2$.
In sequel of this paper, we are going to discuss these cases separately.

\subsection{The settlement of the case $\mathfrak{C}_1$.}

\begin{lem}\label{lemma4.5}   If the Case $\mathfrak{C}_1$ occurs, then $M^3$ is Calabi affine equivalent to an open part of the hypersurfaces $Q(c_1,\cdots,c_r;3),\; 1\leq r\leq 3$.
\end{lem}
\textbf{Proof.} Denote
$$A_{ij}^{k}:=A(e_i,e_j,e_k)\equiv A_{ijk},$$
and put
$$a:=A_{222},\qquad b:=A_{233},\qquad c:=A_{333},\qquad d:=A_{223}.$$


By $\mu_2=\mu_3$, we can further choose $e_2$ as a unit vector for which the function $F$, restricted to $\{e_1^\perp\}\bigcap U_pM^3$, attains its maximum $A_{222}\geq0$. It follows that $A_{223}=0$, and $A_{222}\geq2A_{233}$. Thus we get
\begin{align}
(A^1_{ij})={\left( \begin{array}{ccc}
\mu_1&0&0\\
0&0&0\\
0&0&0
\end{array}
\right)},\;\;
(A^2_{ij})={\left( \begin{array}{ccc}
0&0&0\\
0&a&0\\
0&0&b
\end{array}
\right)},\;\;
(A^3_{ij})={\left( \begin{array}{ccc}
0&0&0\\
0&0&b\\
0&b&c \tag{4.5}
\end{array}
\right)}.
\end{align}
By a direct calculation, we have
\begin{align}
R_{22}&=R_{33}=b(b-a),\notag\\
R_{11}&=R_{12}=R_{13}=R_{23}=0.\tag{4.6}
\end{align}
By (\ref{2.11}) and (\ref{2.12}), it yields
$$0=R_{22}(T_2)^2+R_{33}(T_3)^2=
\frac{1}{9}b(b-a)[(a+b)^2+c^2],\eqno(4.7)$$
\begin{align}
0=&\sum (R_{ijkl})^2+R_{22}\sum (A_{2pq})^2+R_{33}\sum(A_{3pq})^2\notag\\
 =&\sum (R_{ijkl})^2+b(b-a)(a^2+3b^2+c^2).\tag{4.8}
\end{align}
If $b>0$, it contradicts to (4.7). If $b<0$, it also contradicts to (4.8). Therefore $b=0$. By (4.8) we get $R_{ijkl}(p)=0$. Since the arbitrary of point $p$, we have the Calabi metric is flat. Combining $\nabla A=0$ and Theorem 1.3, one can get the following classification results:
\begin{enumerate}\item  if $a=0$, $c=0$, then $ M^3$ is Calabi affine equivalent to an open part of the hypersurface $ Q(c_1;3)$;\\
  \item   if $ a\neq0$, $c=0$,   then $ M^3$ is Calabi affine equivalent to an open part of the hypersurface $Q(c_1,c_2;3)$;\\
  \item   if $a\neq0$, $c\neq0$,  then $M^3$ is Calabi affine equivalent to an open part of the hypersurface $Q(c_1,c_2,c_3;3)$.
\end{enumerate}
\hfill $\Box$
\section{classification of case $\mathfrak{C}_2$}
By $\mu_1=2\mu_2>0$, we know $A_{222}=0$. Thus we have
\begin{align}
(A^1_{ij})={\left( \begin{array}{ccc}
\mu_1&0&0\\
0&\mu_2&0\\
0&0&0
\end{array}
\right)},
(A^2_{ij})={\left( \begin{array}{ccc}
0&\mu_2&0\\
\mu_2&0&d\\
0&d&b
\end{array}
\right)},
(A^3_{ij})={\left( \begin{array}{ccc}
0&0&0\\
0&d&b\\
0&b&c
\end{array}
\right)}.
\end{align}
By (\ref{2.9}), we obtain
\begin{align}
R_{11}=&-\mu_2^2,\quad &R_{22}&=-\mu_2^2+b^2+d^2-cd,\quad &R_{33}=&b^2+d^2-cd,\notag\\
R_{12}=&-\mu_2b,\quad &R_{13}&=\mu_2d,\quad &R_{23}=&0.\notag
\end{align}
Using $\nabla A=0$ and the Ricci identity, we have
\begin{align}\label{5.2} 0=&\;A_{223,13}-A_{223,31}=2A_{p23}R_{p213}
+A_{22p}R_{p313}=2b^2\mu_2. \\
0=&\;A_{222,12}-A_{222,21}=3A_{22p}R_{p212}=3\mu_2(d^2-\mu_2^2).
 \\
0=&\;A_{123,23}-A_{123,32}=A_{p23}R_{p123}+
A_{1p3}R_{p223}+A_{12p}R_{p323} \\
 =&\;\mu_2(2b^2+2d^2-cd).\notag
\end{align}
By (\ref{5.2}), (5.3) and (5.4), we obtain
$$b=0,\qquad d^2=\mu_2^2\neq0,\qquad c=2d.$$
Thus the Pick invariant and the scalar curvature are  \be \label{5.5} J=\frac{7}{3}\mu_2^2,\qquad R=-4\mu_2^2.\ee

Now put tangent vectors
\be\tilde e_1:=\frac{\sqrt 2}{2}(e_1+e_3),\quad \tilde e_3 :=\frac{\sqrt 2}{2}(-e_1+e_3),\ee
then $\{\tilde e_1 , e_2 , \tilde e_3\}$ forms an orthonormal basis of $T_pM^3$,  with respect to which, the Fubini-Pick  tensor $A$ takes the following form:
\be A(\tilde e_1,\tilde e_1)= \sqrt 2 \mu_2 \tilde e_1;\;  A(\tilde e_1,e_2)=\sqrt 2 \mu_2 e_2;\; A(\tilde e_1,\tilde e_3)= \sqrt 2 \mu_2 \tilde e_3, \ee
and
$$A(e_2, e_2)=\sqrt 2 \mu_2 \tilde e_1;\;\; A( e_2,\tilde e_3)=0; \;\; A(\tilde e_3,\tilde e_3)=\sqrt 2 \mu_2 \tilde e_1. $$
By parallel translation along geodesics (with respect to the Levi-Civita connection $\nabla$) through $p$ to a normal neighborhood around $p$, we can extend $\{\tilde e_1,e_2 ,\tilde e_3\}$ to obtain a local orthonormal basis $\{E_1,E_2,E_3\}$ on a neighborhood of $p$ such that
\be A(E_1,E_1)= \sqrt 2 \mu_2 E_1;\; \; A(E_1,E_2)=
 \sqrt 2 \mu_2 E_2;\;\; A(E_1,E_3)= \sqrt 2 \mu_2 E_3\ee
holds at every point in a normal neighborhood.
Denote by $\omega^j_i $ the connection form with respect to the orthonormal frame $\{E_i\}$. By $\nabla A=0$,
$$ A_{11i,j}\omega^j=dA_{11i}-2A_{j1i}\omega^j_1
-A_{11j}\omega^j_i,$$
and choose $i=3$, we have \be\label{5.9}\omega^1_3=0.\ee
Similar, by
$$A_{22i,j}\omega^j=dA_{22i}-2A_{j2i}\omega^j_2-A_{22j}
\omega^j_i,$$
and choose $i=2$, we have \be\label{5.10}\omega^1_2=0.\ee
Then (\ref{5.9}) and (\ref{5.10}) show that $E_1$ is a parallel vector field with respect to the Levi-Civita connection.
Thus, by (\ref{5.5}),  we have
\be R_{2323}=\frac{1}{2}R=-2\mu_2^2=const.\ee

By the above these equalities we have the following lemma:

\begin{lem} We have\\
(i) $\nabla E_1= 0$;\\
(ii) $\langle\nabla_{E_i} E_j , E_1 \rangle= 0$, for any $i,j=2,3$.
\end{lem}

This lemma tell us that the distribution by $\mathcal{D}_1:=\{R E_1 \}$ and $\mathcal{D}_2:=span \{ E_2, E_3 \}$ are totally geodesic. Therefore it follows from the de Rham decomposition theorem (\cite{KN}, pp.187) that as a Riemannian manifold, $(M^3,G)$ is locally isometric to a Riemannian product $R \times H^2(-2\mu_2^2)$, where $H^2(-2\mu_2^2)$ is the hyperbolic  plane of constant negative curvature $-2\mu_2^2$, and after identification, the local vector field $E_1$ is tangent to $R$ and $\mathcal{D}_2$ is tangent to $H^2(-2\mu_2^2)$.

Denote by $x=(x_1,x_2,x_3,x_4)^t$ the position vector of $M^3$ in $A^4$. Using the standard parametrization of the hypersphere model of $H^2(-2\mu_2^2)$, we see that there exists local coordinates $(y_1,y_2,y_3 )$ on $M^3$, such that the metric is given by
\be G=(dy_1)^2+(dy_2)^2+\sinh^2( \sqrt{2}\mu_2 y_2)(dy_3)^2, \ee
and $E_1=\frac{\partial x}{\partial y_1}$, and $\frac{\partial x}{\partial y_2}$, $ (\sinh(\sqrt{2}\mu_2 y_2))^{-1}\frac{\partial x}{\partial y_3}$, form a G-orthonormal basis. We may assume that $E_2=\frac{\partial x}{\partial y_2}$
and $(\sinh(\sqrt{2}\mu_2 y_2))E_3=\frac{\partial x}{\partial y_3}.$
Then a straightforward computation shows that
\begin{align}  \nabla _{\frac{\partial x}{\partial y_2}}\frac{\partial x}{\partial y_2}=\;&0,\\
 \nabla _{\frac{\partial x}{\partial y_2}}\frac{\partial x}{\partial y_3}=\;& \nabla _{\frac{\partial x}{\partial y_3}}\frac{\partial x}{\partial y_2}=  \sqrt{2}\mu_2\coth ( \sqrt{2}\mu_2 y_2)\frac{\partial x}{\partial y_3}, \\
 \nabla _{\frac{\partial x}{\partial y_3}}\frac{\partial x}{\partial y_3}=\;&-\sqrt{2}\mu_2 \sinh ( \sqrt{2}\mu_2 y_2)\cosh ( \sqrt{2}\mu_2 y_2) \frac{\partial x}{\partial y_2} .\end{align}
Using the definition of $A$, we get the following system of differential equations, where, in order to simplify the equations, we have put $c=\sqrt{2}\mu_2$ and $Y=(0,0,0,1)^t$.
\begin{align}\frac{\partial^2 x}{\partial y_1\partial y_1}=\; &c \frac{\partial  x}{\partial y_1 } +Y,\\
 \frac{\partial^2 x}{\partial y_1\partial y_2}=\; &c \frac{\partial  x}{\partial y_2 },  \\
\frac{\partial^2 x}{\partial y_1\partial y_3}=\; &c \frac{\partial  x}{\partial y_3 },  \\
\frac{\partial^2 x}{\partial y_2\partial y_2}=\; &c \frac{\partial  x}{\partial y_1 } +Y,\\
\frac{\partial^2 x}{\partial y_2\partial y_3}=\;&c\coth(cy_2) \frac{\partial  x}{\partial y_3 },\\
\frac{\partial^2 x}{\partial y_3\partial y_3}=\;& c\sinh^2(cy_2) \frac{\partial x}{\partial y_1 }-c\sinh(cy_2)\cosh(cy_2)\frac{\partial x}{\partial y_2} +\sinh^2(cy_2) Y.
\end{align}
To solve the above equations, first we solve its corresponding system of homogeneous equations.
\begin{align}\frac{\partial^2 x}{\partial y_1\partial y_1}=\; &c \frac{\partial  x}{\partial y_1 } ,\\
 \frac{\partial^2 x}{\partial y_1\partial y_2}=\; &c \frac{\partial  x}{\partial y_2 },  \\
\frac{\partial^2 x}{\partial y_1\partial y_3}=\; &c \frac{\partial  x}{\partial y_3 },\\
\frac{\partial^2 x}{\partial y_2\partial y_2}=\; &c \frac{\partial  x}{\partial y_1 },\\
\frac{\partial^2 x}{\partial y_2\partial y_3}=\;&c\coth(cy_2) \frac{\partial  x}{\partial y_3 },\\
\frac{\partial^2 x}{\partial y_3\partial y_3}=\;& c\sinh^2(cy_2) \frac{\partial x}{\partial y_1 }-c\sinh(cy_2)\cosh(cy_2)\frac{\partial x}{\partial y_2} .
\end{align}
From (5.22), we know that there exist vector valued functions $P_1(y_2,y_3)$ and $P_2(y_2,y_3)$ such that
\be
x=P_1(y_2,y_3)e^{cy_1}+P_2(y_2,y_3).\ee
From (5.23) and (5.24) it then follows that the vector function $P_2$ is independent of $y_2$ and $y_3$. Hence there exists a constant vector $A_1$ such that $P_2(y_2,y_3)=A_1$. Next, it follows from (5.25) that $P_1(y_2,y_3)$ satisfies that the following differential equation:
\be \frac{\partial^2 P_1}{\partial y_2\partial y_2 }=c^2 P_1.\ee
Hence we can write
\be P_1(y_2,y_3)=Q_1(y_3)\cosh(cy_2)+Q_2(y_3)\sinh(cy_2).\ee
From (5.26), we then deduce that there exists a constant vector $A_2$ such that $Q_1(y_3)=A_2$. The last formula (5.27) implies there exist constant vectors  $A_3$ and $A_4$ such that
\be
Q_2(y_3)=A_3\cos(cy_3)+A_4\sin(cy_3).
\ee
Therefore the general solution of system (5.22-5.27) are
\be x=e^{cy_1}(A_2\cosh(cy_2)+[A_3\cos(cy_3)+A_4\sin(cy_3)] \sinh(cy_2))+ A_1,\ee
where $A_i$ are constant vectors.
On the other hand, we know that $$\bar x =\left(0,0,0,-\frac{y_1}{c}\right)^t$$ is a special solution of equations (5.16-5.21).
Therefore  the  general solutions of equations (5.16-5.21) are
\be x=e^{cy_1}\left\{A_2\cosh(cy_2)+[A_3\cos(cy_3)+A_4\sin(cy_3)] \sinh(cy_2)\right\}+ A_1+\bar x .\ee

Since $M^3$ is nondegenerate, $x-A_1$ lies linearly full in $A^4$.
Hence $A_2,A_3,A_4$ and $(0,0,0,1)$ are linearly independent vectors. Thus there exists an affine transformation $\phi\in SA(4)$ such that
$$A_1=(0,0,0,0)^t,\, A_2=(1,0,0,0)^t,\,A_3=(0,1,0,0)^t,\,A_4=(0,0,1,0)^t.$$
Then the position vector \be x=\left(\cosh(cy_2)e^{cy_1},\;\cos(cy_3)\sinh(cy_2)e^{cy_1},\;
\sin(cy_3)\sinh(cy_2)e^{cy_1},\;-\frac{y_1}{c}\right)^t.\ee
It follows that, up to an affine transformation $\phi\in SA(4)$, $M^3$ locally lies on the graph hypersurface of function
\be x_4=-\frac{1}{2c^2}\ln(x_1^2-(x_2^2+x_3^2)).
\ee
Thus we finally arrive at the following lemma.

\begin{lem}\label{lemma5.2}   If the Case $\mathfrak{C}_2$ occurs, then $M^3$ is Calabi affine equivalent to an open part of the hypersurface $$x_4=-\frac{1}{2c^2}\ln(x_1^2-(x_2^2+x_3^2)),$$
where the constant $-2c^2$ is the scalar curvature of $M^3$.
\end{lem}

Combining Lemma 4.3, Lemma 4.4, Lemma 4.5 and Lemma 5.2, we complete the proof of Theorem 1.5. \hfill $\Box$


\end{document}